\documentclass[reqno]{amsart}
\usepackage[utf8]{inputenc}
\usepackage{amssymb, amsmath, amstext, amsopn, amsthm, amscd}

\usepackage[numbers]{natbib}

\usepackage{tikz}
\usepackage{graphicx}
\usetikzlibrary{positioning}
\usepackage{tikz-cd}
\usepackage{algpseudocode}
\usepackage{algorithm}

\newcommand{\Cc}{\mathbb{C}}
\newcommand{\Rr}{\mathbb{R}}

\newcommand{\Ss}{\mathbb{S}}
\newcommand{\CC}{\mathfrak{c}}
\newcommand{\GG}{\mathfrak{g}}
\newcommand{\NN}{\mathfrak{n}}
\newcommand{\GL}{\mathfrak{gl}}
\newcommand{\SL}{\mathfrak{sl}}

\newcommand{\SO}{\mathfrak{so}}
\newcommand{\SP}{\mathfrak{sp}}

\newcommand{\ZZ}{\mathfrak{z}}

\DeclareMathOperator{\Tr}{Tr}

\DeclareMathOperator{\Sym}{Sym}
\newcommand{\trans}{\dagger}

\theoremstyle{plain} 
\newtheorem{thm}{Theorem}[]
\newtheorem{cor}{Corollary}[]
\newtheorem{lem}{Lemma}[]

\theoremstyle{definition}
\newtheorem{defn}{Definition}[]

\newtheorem{rem}{Remark}[]

\title{A minimal-variable symplectic method for isospectral flows}
\author{Milo Viviani}
\begin{document}






\maketitle

\begin{abstract}
Isospectral flows are abundant in mathematical physics; the rigid body, the the Toda lattice, the Brockett flow, the Heisenberg spin chain, and point vortex dynamics, to mention but a few.
Their connection on the one hand with integrable systems and, on the other, with Lie--Poisson systems motivates the research for optimal numerical schemes to solve them.
Several works about numerical methods to integrate isospectral flows have produced a large varieties of solutions to this problem.
However, many of these algorithms are not intrinsically defined in the space where the equations take place and/or rely on computationally heavy transformations. In the literature, only few examples of numerical methods avoiding these issues are known, for instance, the \textit{spherical midpoint method} on $\SO(3)$.   
In this paper we introduce a new minimal-variable, second order, numerical integrator for isospectral flows intrinsically defined on quadratic Lie algebras and symmetric matrices. The algorithm is isospectral for general isospectral flows and Lie--Poisson preserving when the isospectral flow is Hamiltonian. 
The simplicity of the scheme, together with its structure-preserving properties, makes it a competitive alternative to those already present in literature.
\end{abstract}

\section{Introduction} 
The numerical integration of isospectral flows is a classical subject of study in numerical analysis \cite{hlw,IsMkNoZa2000}. 
The interest in this problem is motivated by the numerical simulation of integrable systems, which are deeply related to isospectral flows via the Lax pair formulation. 
The quasi-periodic dynamics of integrable systems depends on the presence of a large number of first integrals. 
In the Lax pair formulation, some of these first integrals can be presented as a linear combination of the eigenvalues of the dynamical variable. 
Therefore, the preservation of the spectrum of the dynamical variable is a key feature of a numerical scheme applied to isospectral flows, in order to expect the right qualitative behaviour of the discrete approximate solutions \cite{hlw}. 
Furthermore, as a special case, Lie--Poisson systems on the dual of reductive Lie algebras can be seen as isospectral flows \cite{ModViv2019}. 
A reductive Lie algebra is defined as the direct sum of a semisimple Lie algebra and an abelian Lie algebra. 
In this paper, the crucial property of real semisimple Lie algebras is that they can be represented as matrix Lie algbras which are closed under conjugate transpose \cite[Prop. 6.28]{kna}. Moreover, any real matrix Lie algebra which is closed under conjugate transpose is reductive \cite[Prop. 1.56]{kna}.\footnote{In view of \cite[Sec. I.8]{kna}, there is no restriction in looking at real matrix Lie algebras, since the complex ones can be seen as real matrix Lie algebras of double dimension.}
Because of this, Lie--Poisson systems on the dual of a reductive Lie algebra can be equivalently seen as isospectral flows of the form \eqref{eq:lie_poiss_gen} below.
Lie--Poisson systems originate from the Poisson reduction of canonical Hamiltonian systems on the cotangent bundle of a Lie group \cite{MaRa1999}. 
Classical examples of Lie--Poisson systems are the rigid body~\cite{Po1901}, the heavy top and the incompressible Euler equations~\cite{ArKh1998}.
It is known from the backward error analysis, that Lie--Poisson preserving numerical schemes are superior to standard methods when applied to Lie--Poisson systems, especially for long-time simulations.
As state-of-the-art, well established theories on numerical methods for both isospectral and Lie--Poisson systems exist in the literature (see for example \cite{CaIsZa1997, IsMkNoZa2000, MK1996}). 
For Lie--Poisson systems various symplectic algorithms have been developed (see \cite{Mdd2018} for a recent survey).
However, few examples of numerical schemes that are intrinsically defined in the space where the dynamics takes place are known (e.g. \cite{McMoVe2014d}). This issue often causes a lack of efficiency for these schemes, which rely on group to algebra maps (e.g. the matrix exponential or the Cayley map) or a large number of unknowns.
Before presenting our result, let us introduce the mathematical setup used throughout the paper. 

Isospectral flows are first order ODEs of the form:
\begin{equation}\label{iso_gen}
\begin{array}{ll}
	&\dot W = [B(W),W], \quad W\in V\subset\GL(n,\Cc)\\
	&W(0) = W_0.
	\end{array}
\end{equation}
Here, $[\cdot,\cdot]$ denotes the matrix commutator, $V$ is a linear subspace of the Lie algebra $\GL(n,\Cc)$, and the function $B\colon V\to \NN(V)$ maps $V$ into its $\GL(n,\Cc)$-normalizer algebra $\NN(V)$\footnote{$\NN(V)$ is the largest Lie subalgebra of $\GL(n,\Cc)$ such that $[\NN(V),V]\subset V$ (see Definition~\ref{def_normalizer} in Section~2).}.
The most studied case in literature is when $V=\Sym(n,\Rr)$ is the space of symmetric real matrices, for which the normalizer is the Lie algebra of skew-symmetric real matrices $\NN(V) = \SO(n)$.
For Lie--Poisson systems on the dual of a reductive Lie algebra, we have that $V=\GG^*$ is the dual of a reductive Lie subalgebra of $\GL(n,\Cc)$, for which the $\GL(n,\Cc)-$normalizer is $\NN(V) = \GG_0\oplus \CC(\GG)$ (see Definition~\ref{def_centralizer} and Lemma~\ref{lemma_normalizer}, in Section~2).
Throughout the paper, we identify $\GG^*$ with $\GG$, via the Frobenius inner product $\langle A,B\rangle=\Tr(A^\trans B)$, where $^\trans$ is the conjugate transpose. Via this identifications, Lie--Poisson systems on the dual of a reductive Lie algebra $\GG$ take the form:
\begin{equation}\label{eq:lie_poiss_gen}
\begin{array}{ll}
	&\dot W = [\nabla H(W)^\trans,W], \quad W\in \GG\subset\GL(n,\Cc)\\
	&W(0) = W_0,
	\end{array}
\end{equation}
for $H:\GG\rightarrow \Cc$, a smooth function called \textit{Hamiltonian}.

A class of numerical methods to solve \eqref{iso_gen}-\eqref{eq:lie_poiss_gen}, called \textit{Isospectral Symplectic Runge-Kutta} (IsoSyRK), has been introduced in \cite{ModViv2019}. In the case of Lie--Poisson systems these schemes are symplectic. 
In this paper, we focus on the IsoSyRK associated to the implicit midpoint method, which turns out to have a specially nice structure. On the one hand, we provide a simpler proof (avoiding the use of the B-series theory) that for the implicit midpoint method, the respective IsoSyRK defined in \cite{ModViv2019} is isospectral for any $B=B(W)$. On the other hand, we derive a simpler scheme reducing the number of unknowns up to minimality, revealing an intrinsic relation between the implicit midpoint method and the Cayley transform. The resulting integrator, although implicit, is second order, isospectral and symplectic when the isospectral flow is Lie--Poisson. 
The scheme is also intrinsically defined on $V$, for a large class of isospectral flows (see section 2 for details) and, when the isospectral flow is Lie--Poisson, it preserves the coadjoint orbits. Furthermore, only one evaluation of $B(\cdot)$ and two matrix multiplications per iteration are required, making the scheme very efficient.
In the last section of this paper, we show some numerical examples of our scheme and we compare it with the \textit{spherical midpoint method}, which is another minimal-variable Lie--Poisson integrator on $\Rr^3$. Finally, we show how our scheme looks on $\SL(2,\Rr)$, defining what we call the \textit{hyperbolic midpoint method}.

\medskip

\noindent\textbf{Acknowledgements.} 
The author was supported by EU Horizon 2020 grant No 691070, by the Swedish Foundation for International Cooperation in Research and Higher Eduction (STINT) grant No PT2014-5823, by the Swedish Foundation for Strategic Research grant ICA12-0052, and by the
Swedish Research Council (VR) grant No 2017-05040.
The author would like to thank Klas Modin for the support and the enlightening discussions during the work on this paper.


\section{Main result}\label{sec:main_result}
Let us consider an isospectral flow of the form \eqref{iso_gen}. In order to present our result, we need a short detour on some concepts and basic results on Lie algebras. As already mentioned in section~1, for \eqref{iso_gen} to be well defined, we have to require $B(\cdot)$ to take values in the $\GL(n,\Cc)$-normalizer algebra of $V$. We recall here the definition of the normalizer Lie group and normalizer Lie algebra:
\begin{defn}\label{def_normalizer}
Let $G$ be a Lie group and $\GG$ its Lie algebra. 
Furthermore, let $V\subseteq \GG$ be a linear subspace. 
Then the two sets
\begin{align*}
&N(V)=\lbrace g\in G \mid g^{-1} V g\subseteq V \rbrace\\
&\NN(V)=\lbrace \xi\in \GG \mid [\xi,V]\subseteq V \rbrace
\end{align*}
are respectively called the $G$-normalizer and the $\GG$-normalizer of $V$. 
Notice that $N(V)$ is a subgroup of $G$ and $\NN(V)$ is a Lie subalgebra of $\GG$.
\end{defn}
A related concept to normalizer is the centralizer Lie algebra.
\begin{defn}\label{def_centralizer}
Let $\GG$ be a Lie algebra and let $V\subseteq \GG$ be a linear subspace. 
Then the set
\begin{align*}
&\CC(V)=\lbrace \xi\in \GG \mid [\xi,V]=0 \rbrace
\end{align*}
is called the $\GG$-centralizer of $V$. 
Notice that $\CC(V)$ is a Lie subalgebra of $\GG$.
\end{defn}
We now recall the definition of a $J-$quadratic Lie algebra.
\begin{defn}
A Lie subalgebra $\GG$ of $\GL(n,\Cc)$ is called $J-$quadratic Lie algebra if there exists an invertible matrix $J$, such that
		\begin{equation}\label{eq:quadratic_alg1}
			W\in \GG \iff W^\dagger J + J W = 0.
		\end{equation}
\end{defn}

\begin{lem}
Let $\GG$ be a Lie subalgebra of $\GL(n,\Cc)$ such that there exists a matrix $J$ for which 
		\begin{equation*}
			W\in \GG \iff W^\dagger J + J W = 0.
		\end{equation*}  
Then $\GG=\GG^\trans$ implies $J^2 \in \CC(\GG) $. Moreover, if $J$ is invertible and $J^2 \in \CC(\GG) $, then $\GG=\GG^\trans$.
\end{lem}
\proof
Suppose $\GG=\GG^\trans$. Then, for all $W\in\GG$, both the following identities hold: 
\begin{align*}
&W^\dagger J + J W = 0,\\
&W J + J W^\trans = 0.
\end{align*}
The second of these implies that $W J^2 + J W^\trans J = 0$ and the first one $J W^\dagger J + J^2 W = 0$. Subtracting these identities, we get $[W,J^2]=0$. Hence $J^2 \in \CC(\GG) $.

Now assume that $\GG$ is $J-$quadratic and $J^2 \in \CC(\GG) $. Then, for all $W\in\GG$ we have $0 = JW^\dagger J + J^2 W = JW^\dagger J +  WJ^2$ and hence $W J + J W^\trans = 0$, being $J$ invertible. Therefore $\GG$ and $\GG^\trans$ are defined by the same identity and they coincide. \endproof

\begin{lem}\label{lemma_normalizer}
Let $\GG$ be a Lie subalgebra of $\GL(n,\Cc)$ such that $\GG=\GG^\trans$. Then the \linebreak$\GL(n,\Cc)-$normalizer of $\GG$ is $\NN(\GG)=\GG_0\oplus\CC(\GG)$, where $\GG_0$ is the semisimple ideal of $\GG$ such that $\GG =\GG_0 \oplus \ZZ(\GG)$, for $\ZZ(\GG)$ the center of $\GG$.\footnote{Such a decomposition always exists being $\GG$ reductive by \cite[Prop. 1.56]{kna}.}
\end{lem}
\proof
Let $\GG^\perp$ be the orthogonal complement of $\GG$ in $\GL(n,\Cc)$ with respect to the Frobenius inner product. It is not hard to check that the following properties hold:
\begin{align*}
[\GG,\GG]\subset\GG,\hspace{1cm} 
[\GG,\GG^\perp]\subset\GG^\perp,\hspace{1cm} 
[\GG^\perp,\GG^\perp]\subset\GG.
\end{align*}
Hence, if $A\in \NN(\GG)\cap \GG^\perp$ it must be $[\GG,A]=0$. Therefore, $A$ has to be in $\CC(\GG)$. Moreover, we have that the following inclusions always hold $\ZZ(\GG)\subset\CC(\GG)\subset\NN(\GG)$. Therefore, $\NN(\GG)=\GG_0\oplus\CC(\GG)$, being $\GG_0$ centerless.
\endproof
Notice that, since always $\NN(\GG)^\trans = \NN(\GG^\perp)$, we have that $\NN(\GG^\perp)=\GG_0\oplus\CC(\GG)^\trans$, whenever $\GG=\GG^\trans$. In conclusion, we have the following corollary:
\begin{cor}\label{cor1}
Let $\GG$ be a $J-$quadratic Lie subalgebra of $\GL(n,\Cc)$ such that $J^2 \in \CC(\GG)$. Then the $\GL(n,\Cc)-$normalizer of $\GG$ and $\GG^\perp$ are respectively $\NN(\GG)=\GG_0\oplus\CC(\GG)$ and $\NN(\GG^\perp)=\GG_0\oplus\CC(\GG)^\trans$. In particular, under the identification of $\GG^*$ with $\GG$ via the Frobenius inner product, any Lie--Poisson system on $\GG^*$ can be written in the form \eqref{eq:lie_poiss_gen}.
\end{cor}
Notice that by \cite[Prop. 1.56]{kna} any $\GG$ such as in Corollary~\ref{cor1} is a reductive Lie algebra. As mention in the introduction, Lie--Poisson systems on the dual of reductive Lie algebras can be written in the form \eqref{eq:lie_poiss_gen}. In particular, this is true for Lie--Poisson systems on the dual of $\GG\oplus\mathfrak{Z}$, where $\mathfrak{Z}$ is an Abelian Lie algebra. We can now state the main result of this paper.
\begin{thm}\label{thm1}
Let $W_k\in D\subset V$, for a domain $D$ in the linear subspace $V\subset \GL(n,\Cc)$. Assume that the normalizer splits as $\NN(V)=\GG_0\oplus\CC(V)$, for some Lie algebra $\GG_0$, which satisfies
		\begin{equation}\label{eq:quadratic_alg_gen}
			N\in \GG_0 \iff N^\dagger P + P N = 0.
		\end{equation}
for some constant matrix $P$. Furthermore, let $B:D\subset \GL(n,\Cc)\rightarrow \NN(V)$ be continuously differentiable. Then, for some $h> 0$, there exists $\widetilde{W}\in V$ such that the numerical scheme $W_k\rightarrow W_{k+1}$, implicitly defined by:
\begin{equation}\label{IsoMP_min}
\boxed{
\begin{array}{ll}
&W_{k}=(Id - \frac{h}{2}B(\widetilde{W}))\widetilde{W}(Id + \frac{h}{2}B(\widetilde{W}))\\
&W_{k+1}=(Id + \frac{h}{2}B(\widetilde{W}))\widetilde{W}(Id - \frac{h}{2}B(\widetilde{W})),
\end{array}
}
\end{equation}
is a second order isospectral integrator for \eqref{iso_gen}, for any $k\geq 0$.\footnote{Here $Id$ denotes the $n\times n$ identity matrix.}. Moreover, when (\ref{iso_gen}) is a Lie-Poisson system on $\GL(n,\Cc)^*$ or on the dual of some $J-$quadratic Lie algebra $\GG$ such that $J^2 \in \CC(\GG) $ (or even on $\GG\oplus\mathfrak{Z}$, where $\mathfrak{Z}$ is an Abelian Lie algebra), then (\ref{IsoMP_min}) is a Lie-Poisson integrator for (\ref{iso_gen}) which preserves the coadjoint orbits in $\GG^*$.
\end{thm}

\begin{rem}
The main contribution of Theorem~\ref{thm1}, with respect to the results presented in \citep{ModViv2019}, is that the scheme \eqref{IsoMP_min} is a minimal-variable isospectral (Lie--Poisson) integrator. Minimal-variable means here that the only unknown is $\widetilde{W}$, which lives in a vector space of dimension $\dim(V)$. Hereafter, we will refer to the scheme \eqref{IsoMP_min} as the \textit{isospectral minimal midpoint}. Moreover, the proof of the properties of \eqref{IsoMP_min}, unlikely to \cite[Cor. 1]{ModViv2019}, does not require any application of the B-series theory and reveals a deep connection with the Cayley transform (see the proof of Lemma~\ref{lemma2}). The latter is a quite interesting fact because the Cayley transform arises as a necessary consequence of the use of the implicit midpoint scheme and not, as it has always appeared in literature, as a prescribed choice to construct a certain numerical scheme. We also emphasize that the condition \eqref{eq:quadratic_alg_gen} and the ones on $J$ in Theorem~\ref{thm1} to get a Lie--Poisson integrator is slightly more general than the one considered in \cite[Thm. 1-2]{ModViv2019}.  
\end{rem}
We will give the proof of Theorem~\ref{thm1} in some lemmas. 
\begin{lem}\label{lemma1}
Let $B:D\subset\GL(n,\Cc)\rightarrow \GL(n,\Cc)$ continuously differentiable in the domain $D$. Then, for every $Y\in D$, there exist $\overline{h}>0$ such that the equation
\begin{equation}\label{eq:dcay}
Y=(Id - \frac{h}{2}B(X))X(Id + \frac{h}{2}B(X))
\end{equation}
has a solution $X\in\GL(n,\Cc)$ for any $0\leq h<\overline{h}$.
\end{lem}
\proof
In order to get a solution to \eqref{eq:dcay}, we consider the function $F_h(X) := Y + \frac{h}{2}[B(X),X] + \frac{h^2}{4}B(X)XB(X)$, such that \eqref{eq:dcay} is equal to $X = F_h(X)$. In order to determine $\overline{h}$, we consider the initial value problem:
\begin{align*}
\dfrac{d}{dh}X &= \dfrac{\partial F_h(X)}{\partial h} + DF_h(X)\left[\dfrac{d}{dh}X\right] \\
X(0) &= Y.
\end{align*}
Since $\dfrac{\partial F_h(X)}{\partial h}$ is continuous, if we prove that the operator $(Id - D F_h(X))$ is continuous and invertible, then the Peano existence theorem will ensure a solution $X(h)$, for any $h$ in some interval $[0,\overline{h})$, for $\overline{h}>0$. Indeed, $D F_h(X)=h G(D B (X),B(X),X,h)$, where $G$ is polynomial in its variables and $D B (X)$ is continuous by hypothesis. Hence, $D F_h(X)\rightarrow 0$, for $h\rightarrow 0$, therefore there exist some $\overline{h}>0$ such that $(Id - D F_h(X))$ is invertible, for any $0\leq h<\overline{h}$.
\endproof

\endproof

\begin{lem}\label{lemma2}
Let $B:D\subset\GL(n,\Cc)\rightarrow \GL(n,\Cc)$ continuously differentiable in the domain $D$ and let $0\leq h<\overline{h}$ as in lemma~\ref{lemma1}. Then, for every $W_k\in D$, the numerical scheme $W_k\rightarrow W_{k+1}$ is isospectral. Moreover, if $D\subset V$, for $D$ domain in the linear subspace $V\subset \GL(n,\Cc)$, and $B:D\subset\GL(n,\Cc)\rightarrow \NN(V)$, where $\NN(V)=\GG_0\oplus\CC(V)$, for some Lie algebra $\GG_0$ which satisfies \eqref{eq:quadratic_alg_gen}, then $W_{k+1}\in V$. Furthermore, when $V$ is a $J-$quadratic Lie algebra such that $J^2 \in \CC(\GG) $, then $W_{k+1}\in \mathcal{O}_{W_k}\subset V$, where $\mathcal{O}_{W_k}$ is the coadjoint orbit of which $W_k$ belongs.
\end{lem}
\proof
Clearly $W_{k+1} = (Id + \frac{h}{2}B(\widetilde{W}))(Id - \frac{h}{2}B(\widetilde{W}))^{-1}W_k(Id - \frac{h}{2}B(\widetilde{W}))(Id + \frac{h}{2}B(\widetilde{W}))^{-1}$ and hence $W_{k+1}$ and $W_k$ are similar. Furthermore, we notice that $(Id - \frac{h}{2}B(\widetilde{W}))(Id + \frac{h}{2}B(\widetilde{W}))^{-1}=\mbox{Cay}(\frac{h}{2}B(\widetilde{W}))$, where $\mbox{Cay}$ is the Cayley transform. Therefore we have that:
 \begin{equation}\label{IsoMP_cay}
W_{k+1}=\mbox{Cay}\left(\frac{h}{2}B(\widetilde{W})\right)^{-1}W_k\mbox{Cay}\left(\frac{h}{2}B(\widetilde{W})\right).
\end{equation}
Assuming $\NN(V)=\GG_0\oplus\CC(V)$, for some Lie algebra $\GG_0$ which satisfies \eqref{eq:quadratic_alg_gen}, by \cite[Lemma IV.8.7]{hlw}, $\mbox{Cay}\left(\frac{h}{2}B(\widetilde{W})\right)$ is in the normalizer group $N(V)$ of $V$ and therefore $W_{k+1}$ is in $V$ as well.
When $V = \GG^*$ for $\GG$ a $J-$quadratic Lie algebra such that $J^2 \in \CC(\GG) $, the transformation \eqref{IsoMP_cay} coincides with the coadjoint action of $G$ on $\GG^*$, where $G$ is the respective connected component to the identity of a Lie group with Lie algebra $\GG$. Therefore, \eqref{IsoMP_cay} fixes the coadjoint orbits.
\endproof

\begin{rem}
We point out that the equation \eqref{IsoMP_cay} reveals an interesting relation between the Cayley transform and the implicit midpoint method. Indeed, the fact that the Cayley transform appears as a consequence of the reduction of the implicit midpoint from $T^*GL(n,\Cc)$ to $\GL(n,\Cc)^*$ may indicate a deeper, perhaps canonical, relation between symplectic Runge--Kutta methods and the Cayley transform. The former are associated to conservation of quadratic first integrals of ODEs and the latter to transforming quadratic Lie algebras in quadratic Lie groups. 
\end{rem}

\begin{cor}
Let the hypotesis of Lemma~\ref{lemma1} hold. Then, if $\NN(V)=\GG_0\oplus \CC(V)$  and $\GG_0$ is a compact Lie algebra, there exists $\overline{h}>0$ independent from $k$ such that the scheme \eqref{IsoMP_min} has solution.
\end{cor}
\proof
Since $\GG_0$ is a compact Lie algebra, the associate connected Lie group $G$ is compact and therefore the orbits of the action \eqref{IsoMP_cay} are compact. Hence, we can find a minimum $\overline{h}>0$ in Lemma~\ref{lemma1} independent from the iteration $k\geq 0$.
\endproof

\begin{lem}\label{lemma3}
Let $B:D\subset\GL(n,\Cc)\rightarrow \GL(n,\Cc)$ continuously differentiable in the domain $D$ and let $0\leq h<\overline{h}$ as in lemma~\ref{lemma1}. Then, for every $W_k\in D$, the numerical scheme $W_k\rightarrow W_{k+1}$ in \eqref{IsoMP_min} descends from the method defined in \cite[Def. 1]{ModViv2019} associated with the implicit midpoint method. 
In particular, if $B = \nabla H^\trans$ for some function $H:D\subset\GL(n,\Cc)^*\rightarrow\Rr$, then the method is Lie--Poisson in $\GL(n,\Cc)^*$.
\end{lem}

\proof Consider the second order method as defined in \cite[Def. 1]{ModViv2019} associated with the implicit midpoint method:
\begin{equation}\label{IsoMP}
\begin{array}{llll}
&X=-h(W_k + \frac{1}{2}X)B(\widetilde{W})\\
&Y=hB(\widetilde{W})(W_k + \frac{1}{2}Y)\\
&K=\frac{h}{2}B(\widetilde{W})(X + K)\\
&\widetilde{W}=W_k + \frac{1}{2}(X+Y+ K)\\
&W_{k+1}=W_k + h[B(\widetilde{W}),\widetilde{W}],\end{array}
\end{equation}
for $k\geq 0$ with unknowns $X,Y,K\in\GL(n,\Cc)$. It is not hard to check that the following identities hold:
\begin{align*}
&Y+K=hB(\widetilde{W})\widetilde{W}\\
&X=-h\widetilde{W}B(\widetilde{W})+\frac{h^2}{2}B(\widetilde{W})\widetilde{W}B(\widetilde{W}).
\end{align*}
Applying these to (\ref{IsoMP}) we get, after some computations, the scheme \eqref{IsoMP_min}.
In \cite[Thm. 3]{ModViv2019}, it has been proven that when $B(W)=\nabla H(W)^\trans$, for some functions $H:\GL(n,\Cc)^*\rightarrow\Rr$, the method is a Lie-Poisson integrator in $\GL(n,\Cc)^*$. 
The scheme $W_k\rightarrow W_{k+1}$ defined in \eqref{IsoMP_min} coincides with \eqref{IsoMP}, but with the elimination of the intermediate variables $X,Y,K$. 
Therefore, \eqref{IsoMP_min} is a Lie--Poisson integrator.  
\endproof

\proof $[Theorem~\ref{thm1}]$.
The proof simply follows from the lemmas. Lemma~\ref{lemma1} says that the method \eqref{IsoMP_min} has solution for $h$ sufficiently small. Under the assumptions of Lemma~\ref{lemma2} proves the isospectrality of the scheme and its intrinsic restriction $V$, when $\NN(V)=\GG_0\oplus\CC(V)$, for some Lie algebra $\GG_0$ which satisfies \eqref{eq:quadratic_alg_gen}. Finally, in Lemma~\ref{lemma3} is shown that the scheme descends from the isospectral midpoint method defined in \cite{ModViv2019} and therefore it is a second order Lie--Poisson integrator in $\GL(n,\Cc)^*$, when \eqref{iso_gen} is. Putting together Lemma~\ref{lemma2}, Lemma~\ref{lemma3} and Corollary~\ref{cor1}, we have that when a $V$ is the dual of a $J-$quadratic Lie algebra such that $J^2 \in \CC(\GG)$ (possibly plus a commutative Lie algebra) and \eqref{iso_gen} is Lie--Poisson, the scheme \eqref{IsoMP_min} is a Lie--Poisson integrator on $V$, which preserves the coadjoint orbits.
\endproof

 \begin{rem}
 We notice that the isospectral minimal midpoint \eqref{IsoMP_min} is somehow similar to the \textit{modified implicit midpoint rule} introduced in \cite{CaIsZa1997}. However, in their scheme, $\widetilde{W}$ was set to be $\frac{W_{k+1}+W_k}{2}$ which does not hold in general while solving the isospectral minimal midpoint \eqref{IsoMP_min}. In fact, even though the scheme in \cite{CaIsZa1997} is isospectral, it is not symplectic.
 \end{rem} 
  \begin{rem}
The isospectral minimal midpoint \eqref{IsoMP_min} can be derived in a different way, as proposed in \cite{Mdd2018}. The construction there is more general and \eqref{IsoMP_min} can be recovered choosing as a retraction map the Cayley transform instead of the exponential map. This surprising connection opens up a question about a geometrical description of the methods proposed in \cite[Def. 1]{ModViv2019}. Let us consider for any $s=1,2,\dots$ and $a$ $s\times s$ real matrix, a retraction map $\tau_a:\GG^{\oplus s}\rightarrow G^{\times s}$. Then, similarly to \cite{Mdd2018}, for each $i=1,\dots,s$, it is implicitly defined by the differential of the retraction map $d\tau_a$ a discrete map $W_k\mapsto \widetilde{W}_i\ni\GG$. Finally, we define our integrator $\Psi_h:W_k\mapsto W_{k+1}$ as:
\[
W_{k+1} = W_k + h\sum_{i=1}^s b_i[B(\widetilde{W}_i),\widetilde{W}_i],
\]
for some real numbers $b_i$. The question is whether, for any $s$-stage symplectic Runge--Kutta method, there exists a retraction map $\tau_a:\GG^{\oplus s}\rightarrow G^{\times s}$, such that any Lie--Poisson integrator defined in \cite{ModViv2019} can be obtained in this way.
 \end{rem}

\section{Numerical examples} \label{sec:examples}

In this section we present some applications of the isospectral minimal midpoint \eqref{IsoMP_min} on isospectral flows and Lie--Poisson systems found in literature.
We also compare our method in the case of $\SO(3)\cong(\Rr^3,\times)$ with the \textit{spherical midpoint method}, showing that the isospectral minimal midpoint \eqref{IsoMP_min} has the same computational cost. 
Finally, we show explicitly how the isospectral minimal midpoint \eqref{IsoMP_min} looks on $\SL(2,\Rr)$, applying it to the the point vortex equations on the hyperbolic plane.
For the example considered in this section, we plot the variation of the first integrals of \eqref{iso_gen}. As expected, we get  exact conservation (up to round-off errors) of the Casimir functions and, when the flow is Hamiltonian, near conservation of the Hamiltonian.

\subsection{The generalized rigid body} \label{sec:rigidbody}
A classical example among Hamiltonian isospectral systems is the generalized rigid body. 
It represents a class of completely integrable systems on $\SO(n)$, for every $n\geq 1$, \cite{man}. 
The Hamiltonian is given by
\begin{equation}\label{eq:ham_rb}
H(W) = \frac{1}{2}\Tr((\mathcal{I}^{-1}W)^\trans W),\qquad W\in\SO(n),
\end{equation}
where $\mathcal{I}\colon\SO(n)\rightarrow\SO(n)$ is a symmetric positive definite inertia tensor.
The equations of motion are then
\begin{equation}\label{rigid_body}
\left. \begin{array}{ll} 
	\dot{W} = -[\mathcal{I}^{-1}W,W] \\
	 W(0)=W_0.
			\end{array}\right.
\end{equation} 

We discretize this system for $n=10$.
Our implementation uses Newton iterations for the non-linear system.\footnote{Fixed-point iteration could be used as well, no numerical issues have arisen in our experiments with it.}
The inertia tensor is given by 
\begin{equation}
	(\mathcal{I}^{-1}W)_{ij} =\left\lbrace \begin{array}{ll}  \frac{W_{ij}}{i} , \quad i=1,\ldots,5,j=1,\ldots,10\\
	\frac{W_{ij}}{11-i} , \quad i=6,\ldots,10,j=1,\ldots,10\end{array}\right.
\end{equation}
and we use the stepsize $h=0.1$.
The initial conditions are given by
\begin{equation}
	(W_0)_{ij}=1/10 \quad\text{for}\quad i<j \qquad\text{and}\quad W_0^\trans = -W_0
\end{equation}

\begin{figure}[h!]
\begin{minipage}{1\textwidth}
\begin{tikzpicture}
\centering
 \node (img)  {\includegraphics[scale=0.3]{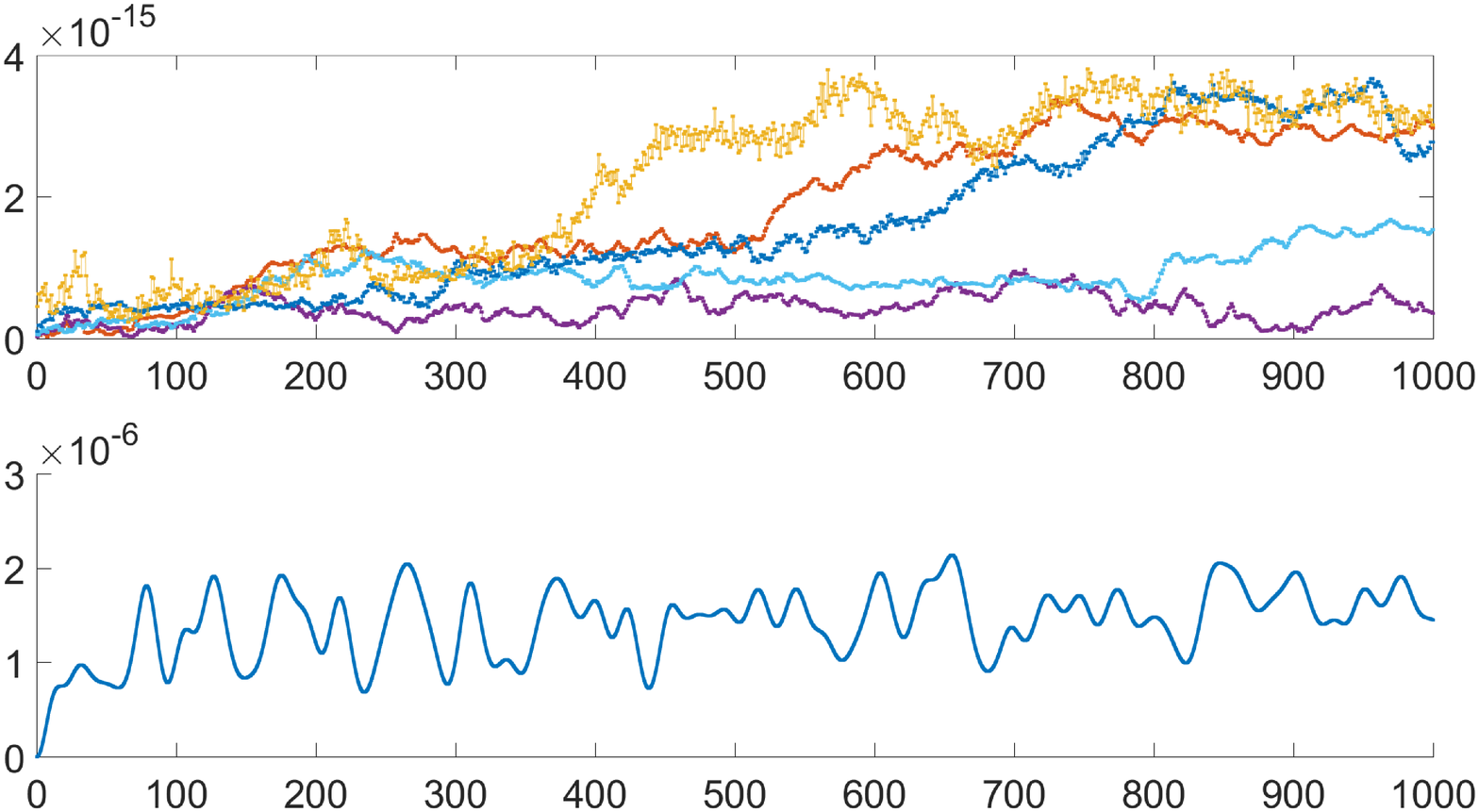}};
  \node[below=of img, node distance=0cm, yshift=4.1cm] {Hamiltonian variation in time};
  \node[below=of img, node distance=0cm, yshift=7.2cm] {Casimir (eigenvalues) variation in time};
 \end{tikzpicture}
\end{minipage}

\caption{Casimir and Hamiltonian variation in time $T=100$, for the generalized rigid body on $\SO(10)$ and time-step $h=0.1$.}\label{fig:rb}
\end{figure}

As shown in Figure~\ref{fig:rb}, the Hamiltonian is nearly conserved and the Casimir functions are conserved up to the accuracy of the Newton iterations. 

\subsection{The Brockett flow}
\begin{figure}[h!]
\begin{minipage}{1\textwidth}
\begin{tikzpicture}
 \node (img)  {\includegraphics[scale=0.3]{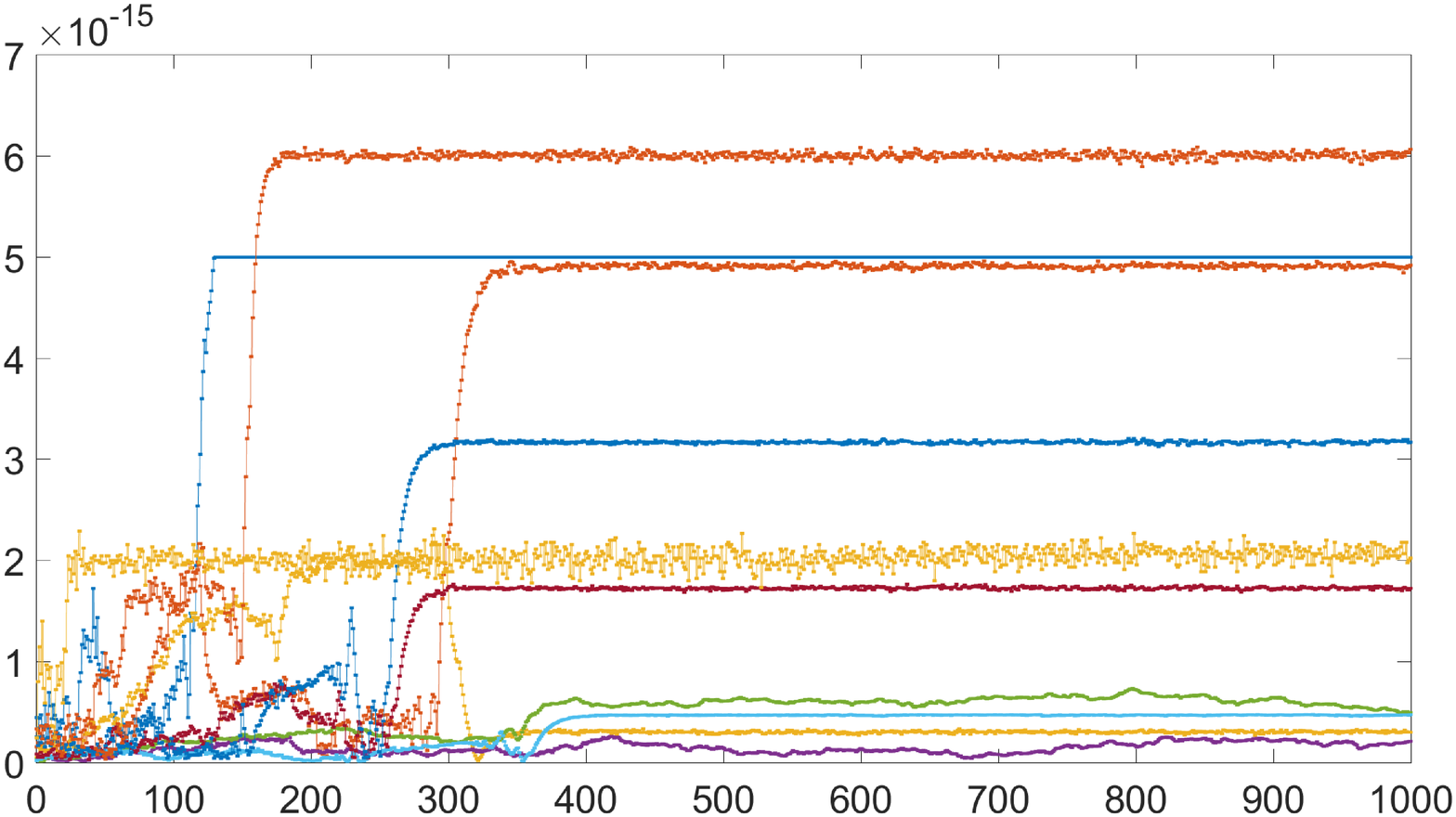}};
\centering
   \node[below=of img, node distance=0cm, yshift=7.4cm] {Casimir (eigenvalues) variation in time};
 \end{tikzpicture}
\end{minipage}
\caption{Eigenvalue variation in time $T=1000$, for the isospectral minimal midpoint \eqref{IsoMP_min} applied to \eqref{eq:Brockett} with time-step $h=0.1$. Then initial value $W_0$ is a randomly generated self-adjoint matrix of dimension $10\times 10$ and $N=diag(1,2,\dots,10)$.}\label{fig:eig_Brockett}
\end{figure}

In this section we specify the isospectral minimal midpoint \eqref{IsoMP_min} for the \textit{Brockett flow}, or \textit{double bracket flow}:
\begin{equation}\label{eq:Brockett}
\dot{W} = [[N,W],W],
\end{equation}
where $N,W$ are $n\times n$ self-adjoint complex matrices.
In \cite{Bro1991}, Brockett shows that for $N$ diagonal matrix with distinct entries and $W_0$ self-adjoint matrix with distinct eigenvalues, for $t\rightarrow\infty$, $W(t)$ converges exponentially fast to a diagonal matrix with the eigenvalues sorted accordingly to the order of the entries of $N$. 
In figure~\ref{fig:eig_Brockett}, we plot the eigenvalue variation for a randomly generated\footnote{Throughout the whole paper, a randomly generated matrix is understood as a matrix whose components are defined as pseudorandom values drawn from the standard uniform distribution on the open interval $(0,1)$, generated via the MATLAB function \texttt{rand}.} self-adjoint initial matrix $W_0$ of dimension $10\times 10$ and $N=diag(1,2,\dots,10)$. The asymptotic stationarity of the eigenvalues variation reflects the fact that $W$ is close to a diagonal matrix.

\subsection{Lie--Poisson systems on $(\Rr^3,\times)$}
On $(\Rr^3,\times)$ the isospectral minimal midpoint \eqref{IsoMP_min} can be written as:
\begin{equation}\label{IsoMP_min_R3}
\begin{array}{ll}
&w_k=\widetilde{w} + \frac{h}{2}\widetilde{w}\times B(\widetilde{w}) - \frac{h^2}{4}B(\widetilde{w})(B(\widetilde{w})\cdot\widetilde{w})\\
&w_{k+1}=\widetilde{w} - \frac{h}{2}\widetilde{w}\times B(\widetilde{w}) - \frac{h^2}{4}B(\widetilde{w})(B(\widetilde{w})\cdot\widetilde{w}),
\end{array}
\end{equation}
for $\widetilde{w},w_k,w_{k+1}\in\Rr^3$ and $B:\Rr^3\rightarrow\Rr^3$. We want to compare the isospectral minimal midpoint \eqref{IsoMP_min_R3} with another minimal-variable symplectic integrator on $\Rr^3$ introduced in \cite{mmv}, i.e. the \textit{spherical midpoint method}:
\begin{equation}\label{SphMP}
\begin{array}{ll}
&w_{k+1}=w_k + h\dfrac{\sqrt{w_{k+1}}\sqrt{w_k}(w_{k+1}+w_k)}{|w_{k+1}+w_k|}\times B\left(\dfrac{\sqrt{w_{k+1}}\sqrt{w_k}(w_{k+1}+w_k)}{|w_{k+1}+w_k|}\right).
\end{array}
\end{equation}
\begin{rem}
Let $B(\cdot)$ be orthogonal with respect to the rays, i.e. $B(w)\cdot w=0$, for every $w\in\Rr^3$. It is immediate to check that then \eqref{IsoMP_min_R3} coincides with the classical midpoint scheme:
\begin{equation}\label{ClassicMP}
\begin{array}{ll}
&w_{k+1}=w_k + h\dfrac{w_{k+1}+w_k}{2}\times B\left(\dfrac{w_{k+1}+w_k}{2}\right).
\end{array}
\end{equation}
In \cite{mmv} it is shown that also \eqref{SphMP} coincides with \eqref{ClassicMP} when $B(\cdot)=\nabla H(\cdot)$, for some Hamiltonian function $H:\Rr^3\rightarrow\Rr$ constant on the rays (which implies $B(\cdot)$ to be orthogonal to the rays). In this case, \eqref{ClassicMP} is known to be symplectic, whereas this fails for general Hamiltonian $H$. Therefore, \eqref{IsoMP_min_R3} can be seen as the second order correction of \eqref{ClassicMP} to be symplectic for any Hamiltonian $H$.
\end{rem}
Let us now consider the two schemes \eqref{IsoMP_min_R3} and \eqref{SphMP}. Both methods are implicit and therefore an implicit solver has to be used. Here we show that they exhibit the same computational cost. The example we consider is the Heisenberg spin chain on $\Rr^{3N}$. For this one has to extend both the isospectral minimal midpoint \eqref{IsoMP_min_R3} and the spherical midpoint \eqref{SphMP} to direct products of $\Rr^3$ (see \cite{mmv} and \cite{ModViv2019}). 

\begin{figure}[htb]
\begin{minipage}{1\textwidth}
\begin{tikzpicture}
 \node (img)  {\includegraphics[scale=0.3]{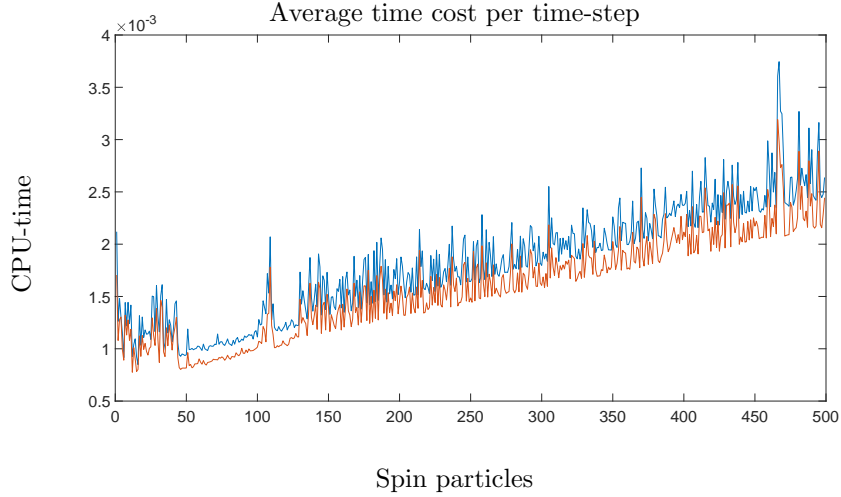}};
\centering
  \node[below=of img, node distance=0cm, yshift=7.2cm] {Average time cost per time-step};
  \node[below=of img, node distance=0cm, yshift=1cm] {Spin particles};
  \node[below=of img, node distance=0cm, xshift=-6cm, yshift=4cm, rotate=90] {CPU-time};
 \end{tikzpicture}
\end{minipage}
\caption{Average time cost per time-step in seconds, with respect to the number of spin particles, for randomly generated initial values. The upper and the lower line are referred to, respectively, the isospectral minimal midpoint \eqref{IsoMP_min_R3} and the spherical midpoint \eqref{SphMP}.}\label{fig:MMPvsSPHMP}
\end{figure}

The Heisenberg spin chain of micromagnetics is defined as:
\begin{equation}\label{HeiSC}
\dot{w}_i = w_i\times(w_{i-1} + w_{i+1}),
\end{equation}
where $w_i\in\Ss^2$, for $i=1,\ldots,N$ and $w_{N+1}=w_1$. It corresponds, up to scaling, to spatial discretization of the Landau-Lifschitz PDE:
\begin{equation}\label{LLPDE}
\dot{w} = w\times\partial_{xx} w,
\end{equation}
for $w:\Ss^1\rightarrow\Ss^2$ smooth. We notice that \eqref{HeiSC} is a Lie--Poisson system on $\Rr^{3N}$, with Hamiltonian:
\[
H(w_1,\ldots,w_N) = \sum_{i=1}^{N} w_i\cdot w_{i+1}.
\]
Clearly, to get a good approximation of \eqref{LLPDE}, $N$ has to be large. 
In figure \ref{fig:MMPvsSPHMP} we show the average time cost for time-step with respect to the number of spin particles for both the isospectral minimal midpoint \eqref{IsoMP_min_R3} and the spherical midpoint \eqref{SphMP}. We conclude that the complexity grows similarly. 

\begin{figure}[h!]
\begin{minipage}{1\textwidth}
\begin{tikzpicture}
 \node (img)  {\includegraphics[scale=0.3]{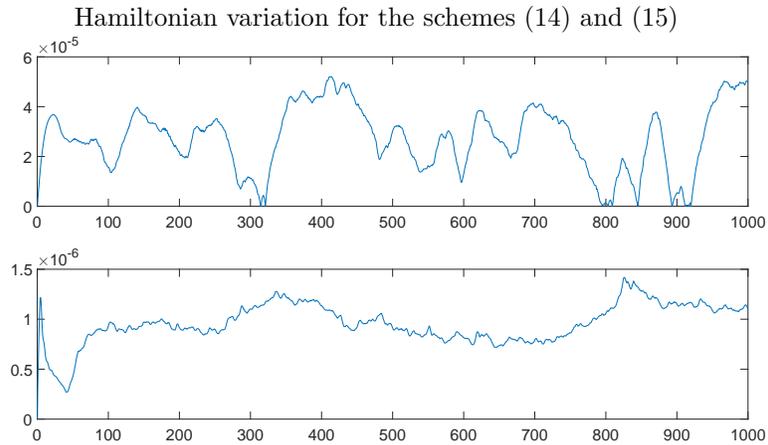}};
\centering
  \node[below=of img, node distance=0cm, yshift=7.4cm] {Hamiltonian variation for the  schemes \eqref{IsoMP_min_R3} and  \eqref{SphMP}};
 \end{tikzpicture}
\end{minipage}
\caption{Hamiltonian variation $|H(k h)-H(0)|$ in time $T=1000$, for the isospectral minimal midpoint \eqref{IsoMP_min_R3} above and the spherical midpoint \eqref{SphMP} below, for $N=100$ spin particles and time-step $h=0.1$.}\label{fig:HM_MMPvsSPHMP}
\end{figure}

In terms of conservation properties, the two schemes exactly preserve the linear invariants and nearly conserve the the quadratic first integrals of the form $\sum_{i,j}w_i^\trans A w_j$. Moreover, the spherical midpoint has the advantage of exactly conserving all the quadratic first integrals of the form $w_i^\trans A w_i$, for some square matrix $A$, whereas the isospectral minimal midpoint \eqref{IsoMP_min_R3} exactly conserves only the quadratic invariants $w_i^\trans w_i$. 
In figure \ref{fig:HM_MMPvsSPHMP} we compare the isospectral minimal midpoint \eqref{IsoMP_min_R3} and the spherical midpoint \eqref{SphMP} for the initial data given as am equispaced discretization in $N=100$ points of the closed spherical curve: 
\[w(x)=(\cos(2\pi x^2)\sin(2\pi x^3),\sin(2\pi x^2)\sin(2\pi x^3),\cos(2\pi x^3)),\]
for $x\in [0,1]$.
We can conclude from figure \ref{fig:HM_MMPvsSPHMP} that the spherical midpoint performs slightly better than the isospectral minimal midpoint \eqref{IsoMP_min_R3}. However, both the schemes show the desired conservation properties due to their symplecticity.
\begin{figure}[h!]
\begin{minipage}{1\textwidth}
\begin{tikzpicture}
 \node (img)  {\includegraphics[scale=0.3]{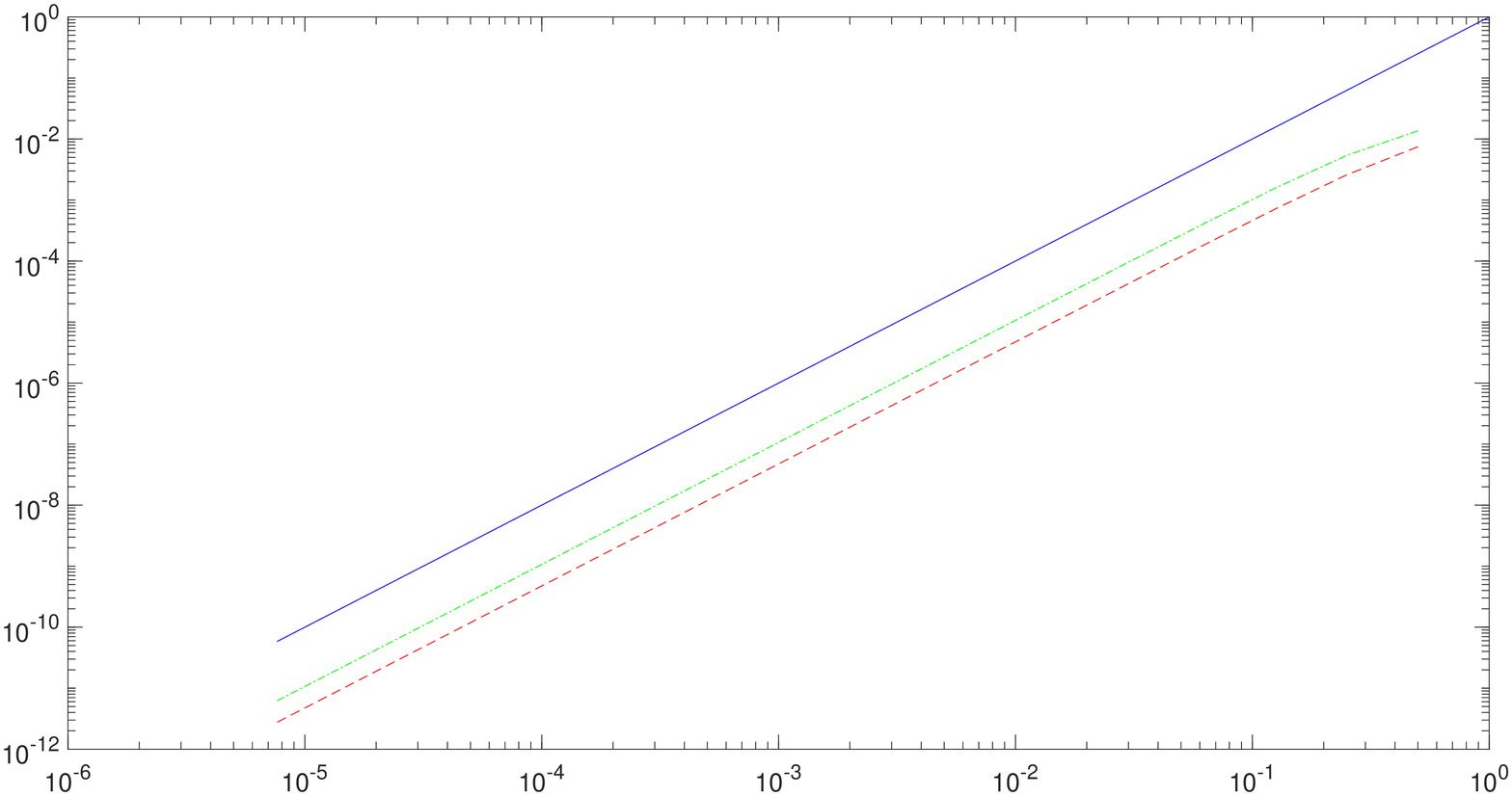}};
\centering
  \node[below=of img, node distance=0cm, yshift=8cm] {Error diagram for the schemes \eqref{IsoMP_min_R3} and  \eqref{SphMP}};
    \node[below=of img, node distance=0cm, yshift=1cm] {Time-step h};
  \node[below=of img, node distance=0cm, xshift=-6cm, yshift=4cm, rotate=90] {Error $\max_k\|w(kh)-w_k\|$};
 \end{tikzpicture}
\end{minipage}
\caption{Maximum error in total time $T=1s$, and time-step $h$, for $h=1, 0.5^{2},\dots,0.5^{17}$, in loglog scale, for the isospectral minimal midpoint \eqref{IsoMP_min_R3}, dot-dashed line, and the spherical midpoint \eqref{SphMP}, dashed line, for $N=100$ spin particles. The continuous line is $h\mapsto h^2$.}\label{fig:err_MMPvsSPHMP}
\end{figure}
Finally, in figure~\ref{fig:err_MMPvsSPHMP}, we present the error diagram for the isospectral minimal midpoint \eqref{IsoMP_min_R3} and the spherical midpoint \eqref{SphMP}. The curves in figure~\ref{fig:err_MMPvsSPHMP} show the expected second order of the schemes, with no significant difference.

\subsection{Lie--Poisson systems on $\SL(2,\Rr)^*$}
In this section we specify the isospectral minimal midpoint \eqref{IsoMP_min} in the case of the Lie algebra $\SL(2,\Rr)$, i.e. the $2\times 2$ matrices with zero trace. The first observation is that $\SL(2,\Rr)\cong \SP(2,\Rr)$, which is a non-compact $J-$quadratic Lie algebra with respect to $J=\begin{bmatrix} 
  0 & -1\\
  1 & 0
\end{bmatrix}.$ From this, it is straightforward to see that, when $B:\GL(2,\Rr)\rightarrow\SL(2,\Rr)$ and $W_k\in\SL(2,\Rr)$, $\widetilde{W}$ in \eqref{IsoMP_min} is also in $\SL(2,\Rr)$. We notice that this is no longer true for $\SL(n,\Rr)$, for $n>2$. On the other hand, any element in $\SL(2,\Rr)$ can be written as a vector in $\Rr^3$, via the vector spaces isomorphism:
\begin{equation}\label{eq:sl2tor3}
\begin{bmatrix} 
  x & y+z\\
  y-z & -x
\end{bmatrix}\mapsto \begin{bmatrix}x\\y\\z
\end{bmatrix},
\end{equation}
for any $x,y,z\in\Rr$. In this coordinates, we can express the isospectral minimal midpoint \eqref{IsoMP_min} for $\SL(2,\Rr)$ as:
\begin{equation}\label{IsoMP_min_sl2}
\begin{array}{ll}
&w_k=\widetilde{w} + \frac{h}{2}2L(\widetilde{w}\times B(\widetilde{w})) - \frac{h^2}{4}M(B(\widetilde{w}),\widetilde{w},B(\widetilde{w}))\\
&w_{k+1}=\widetilde{w} - \frac{h}{2}2L(\widetilde{w}\times B(\widetilde{w})) - \frac{h^2}{4}M(B(\widetilde{w}),\widetilde{w},B(\widetilde{w})),
\end{array}
\end{equation}
for $\widetilde{w},w_k,w_{k+1}\in\Rr^3$ and $B:\Rr^3\rightarrow\Rr^3$, where :
\[
L := \begin{bmatrix} 
  1 & 0 & 0\\
  0 & 1 & 0\\
  0 & 0 & -1
\end{bmatrix}
\]
and
\[
M(a,b,c) :=  \begin{bmatrix}
c_1(a_1b_1 + (b_2 - b_3)(a_2 + a_3)) - (c_2 - c_3)(b_1(a_2 + a_3) - a_1(b_2 + b_3))\\
c_1(b_1(a_2 + a_3) - a_1(b_2 + b_3)) + (c_2 + c_3)(a_1b_1 + (b_2 - b_3)(a_2 + a_3))\\ 
c_1(b_1(a_2 - a_3) - a_1(b_2 - b_3)) + (c_2 - c_3)(a_1b_1 + (b_2 + b_3)(a_2 - a_3))
\end{bmatrix},
\]
for any $a,b,c\in\Rr^3$. 

We notice that the map \eqref{eq:sl2tor3} is a Lie algebra isomorphism from $\SL(2,\Rr)$ to $(\Rr^3,\times_L)$, where $a\times_L b = 2L (a\times b)$, for any $a,b\in\Rr^3$. 
Furthermore, the tensor $L$ defines also the hyperbolic inner product $a\cdot_L b = a\cdot(Lb)$, for any $a,b\in\Rr^3$. 
We recall that the coadjoint orbits in $\SL(2,\Rr)^*$ are hyperboloids of the form $\lbrace x^2+y^2-z^2=const\rbrace$.
Hence, in analogy with the spherical midpoint method in \cite{McMoVe2014d}, we will call \eqref{IsoMP_min_sl2} the \textit{hyperbolic midpoint method}.

We illustrate an application of the hyperbolic midpoint method \eqref{IsoMP_min_sl2} on the \textit{point vortex equations} on the hyperbolic plane (see for example \cite{HwKi2009,HwKi2013,MoNa2014}). The interest in this equations is motivated also by the studies on ideal hydrodynamics on hyperbolic spaces \cite{ChCz2013, KhMi2012}, and in particular on the Euler equations, for which the point vortices can be seen as a finite dimensional approximation \cite{BeMa1982}.
These equations are a Lie--Poisson system on $(\SL(2,\Rr)^*)^N\cong (\Rr^3,\times_L)^N$, with initial values on the coadjoint orbit determined by the equations $w_i\cdot_L w_i=-1$, for $i=1,2,\dots,N$. The Hamiltonian is given by:
\begin{equation}\label{eq:PV_ham} 
H = -\frac{1}{4\pi}\sum_{i\neq j}\Gamma_i\Gamma_j \log\left(\dfrac{w_i\cdot_L w_j + 1}{w_i\cdot_L w_j - 1}\right).
\end{equation}
The equations of motion are then:
\begin{equation}\label{eq:PV_eq}
\dot{x}_i = -\frac{1}{\pi}\sum_{i\neq j}\Gamma_j \dfrac{w_i\times_L w_j}{(w_i\cdot_L w_j)^2 - 1}.
\end{equation}
Equations \eqref{eq:PV_eq} constrain the vortices to move on the hyperboloid $x^2+y^2-z^2=-1$. Furthermore, the $SL(2,\Rr)$ symmetry of \eqref{eq:PV_eq} gives the conservation of the momentum vector:
\begin{equation}\label{eq:PV_mom}
M = \sum_{i=1}^N \Gamma_i w_i.
\end{equation} 
Equations \eqref{eq:PV_eq} and their $SL(2,\Rr)$-relative equilibria have been studied in \cite{HwKi2009, HwKi2013, MoNa2014}. In particular, for two and three vortices most of the stability issues have been worked out in \cite{MoNa2014}. However, unlike to point vortex equations on a sphere \cite{LaMoRo2011}, it is still unknown a general result on the stability of a relative equilibrium of point vortices on the hyperbolic plane. 

\begin{figure}[htb]
\begin{minipage}{1\textwidth}
\begin{tikzpicture}
 \node (img)  {\includegraphics[scale=0.3]{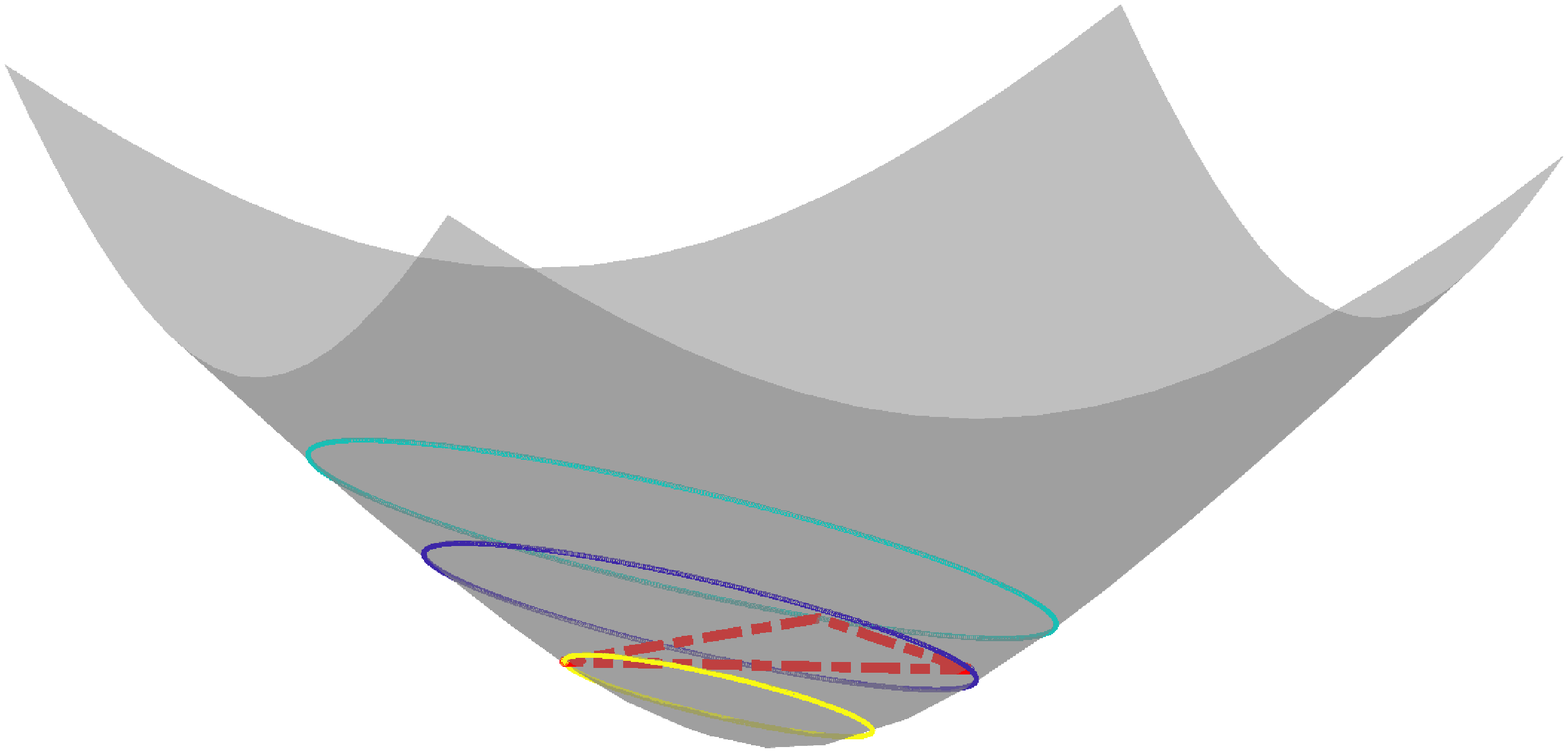}};
\centering
  \node[below=of img, node distance=0cm, yshift=7.2cm] {Point vortex trajectories for $w_0^1$};
 \end{tikzpicture}
\end{minipage}
\begin{minipage}{1\textwidth}
\begin{tikzpicture}
 \node (img)  {\includegraphics[scale=0.3]{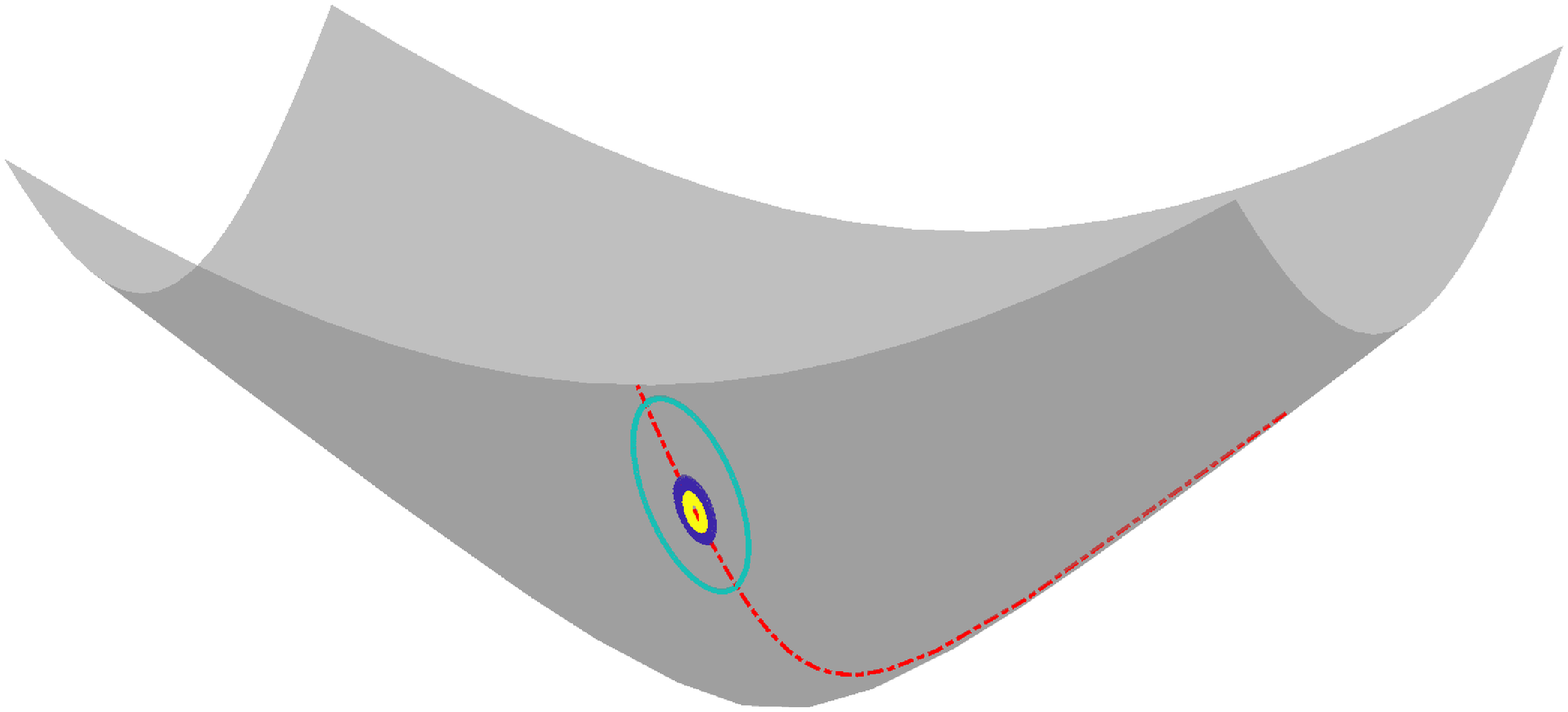}};
\centering
  \node[below=of img, node distance=0cm, yshift=7.2cm] {Point vortex trajectories for $w_0^2$};
 \end{tikzpicture}
\end{minipage}
\caption{Point vortex trajectories for the initial conditions $w_0^1,\Gamma^1$ and $w_0^2,\Gamma^2$, as defined above. The purple, the blue and the yellow lines represent the point vortex trajectories. The red-dashed line is, respectively, the equilateral triangle of the initial configuration $w_0^1$ and the geodesic passing through the initial condition $w_0^2$. The simulations have been carried with the hyperbolic midpoint method \eqref{IsoMP_min_sl2}, with time-step, respectively, $h=0.01$ and  $h=0.001$, total time $T=10$ and $T=1$ and tolerance for the Newton iteration $tol=10^{-13}$.}\label{fig:PV_traj}
\end{figure}

Here we present the results of some numerical simulations of \eqref{eq:PV_eq} with the hyperbolic midpoint method \eqref{IsoMP_min_sl2}. In particular, we consider as initial values some of the relative equilibria found in \cite{HwKi2013,MoNa2014}. Let us take two different initial values in $(\Rr^3)^3$, $w_0^1$ and $w_0^2$, whose columns represent the initial position of three vortices with strengths respectively equal to $\Gamma^1$ and $\Gamma^2$, as defined here below:
 \[
  w^1_0 = \begin{bmatrix}
  -0.5000  & -0.5000 &   1.0000  \\
    0.8660  & -0.8660 &  -0.0000   \\
    1.4142  &  1.4142  &  1.4142 \\
 \end{bmatrix},
 \hspace{1cm}
 w^2_0 = \begin{bmatrix}
  2.6000  &  4.0000  &  3.0000 \\
    0.1923 &   0.1250 &   0.1667 \\
    2.7923  &  4.1250  &  3.1667 \\
 \end{bmatrix},
 \]
 \[
 \Gamma^1 =  \begin{bmatrix} 0.5317 & 0.0761 & 1.0000 \end{bmatrix}, \hspace{2cm} \Gamma^2 =  \begin{bmatrix}0.0990 &   0.8091 & 1.0000\end{bmatrix}.
 \]
The initial condition $w_0^1$ is an equilateral relative equilibrium, whereas  $w_0^2$ is a geodesic relative equilibrium, as defined in \cite{MoNa2014}. The fist initial condition is known to be stable, whereas it is not known for the second one. However, from figure~\ref{fig:PV_traj} we can see that both the initial conditions evolve in close trajectories, which proves numerically the stability for both of them. 

\begin{figure}[htb]
\begin{minipage}{1\textwidth}
\begin{tikzpicture}
 \node (img)  {\includegraphics[scale=0.3]{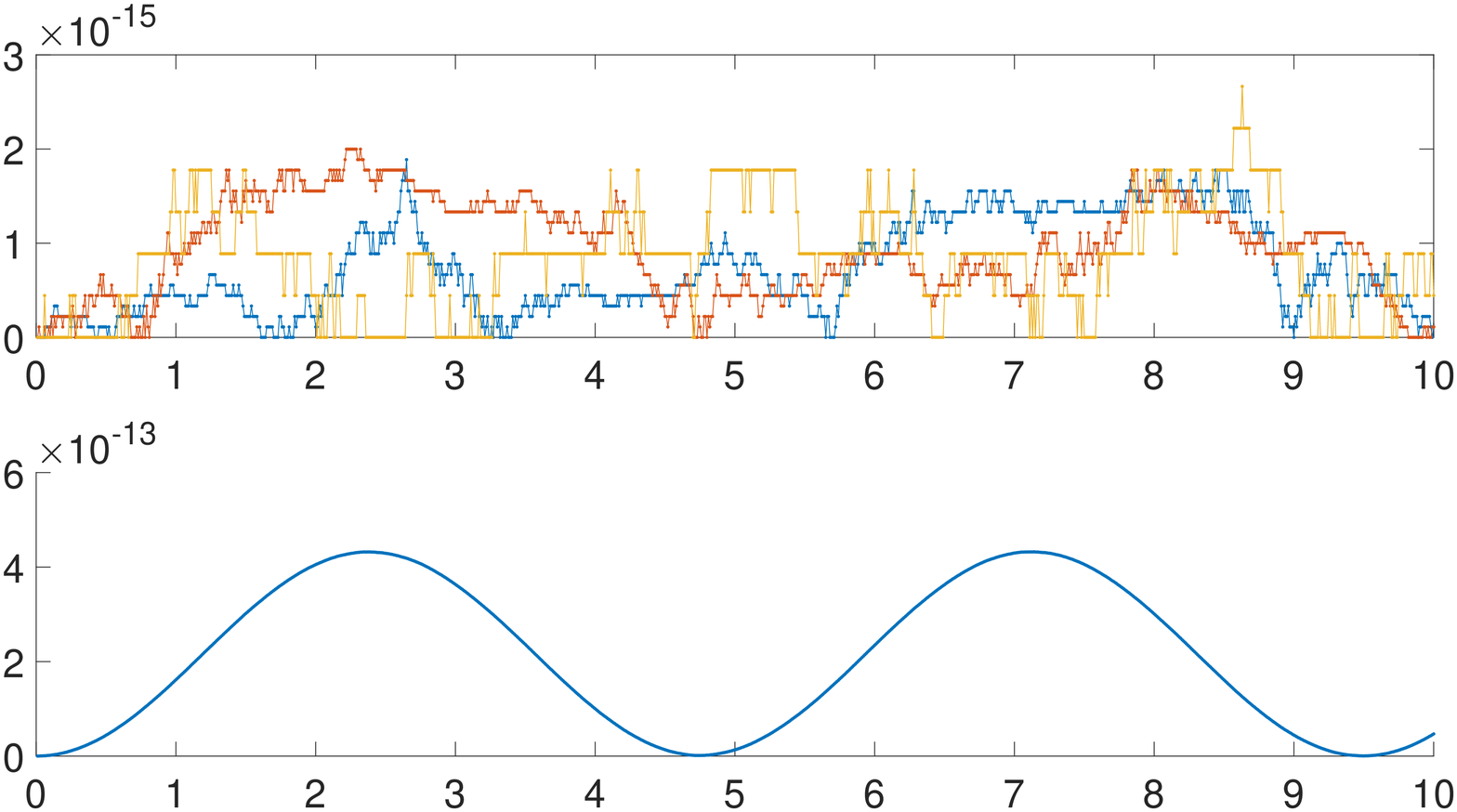}};
\centering
   \node[below=of img, node distance=0cm, yshift=4.1cm] {Hamiltonian variation in time for $w_0^1$};
  \node[below=of img, node distance=0cm, yshift=7.2cm] {Momentum variation in time  for $w_0^1$};
 \end{tikzpicture}
\end{minipage}
\begin{minipage}{1\textwidth}
\begin{tikzpicture}
 \node (img)  {\includegraphics[scale=0.3]{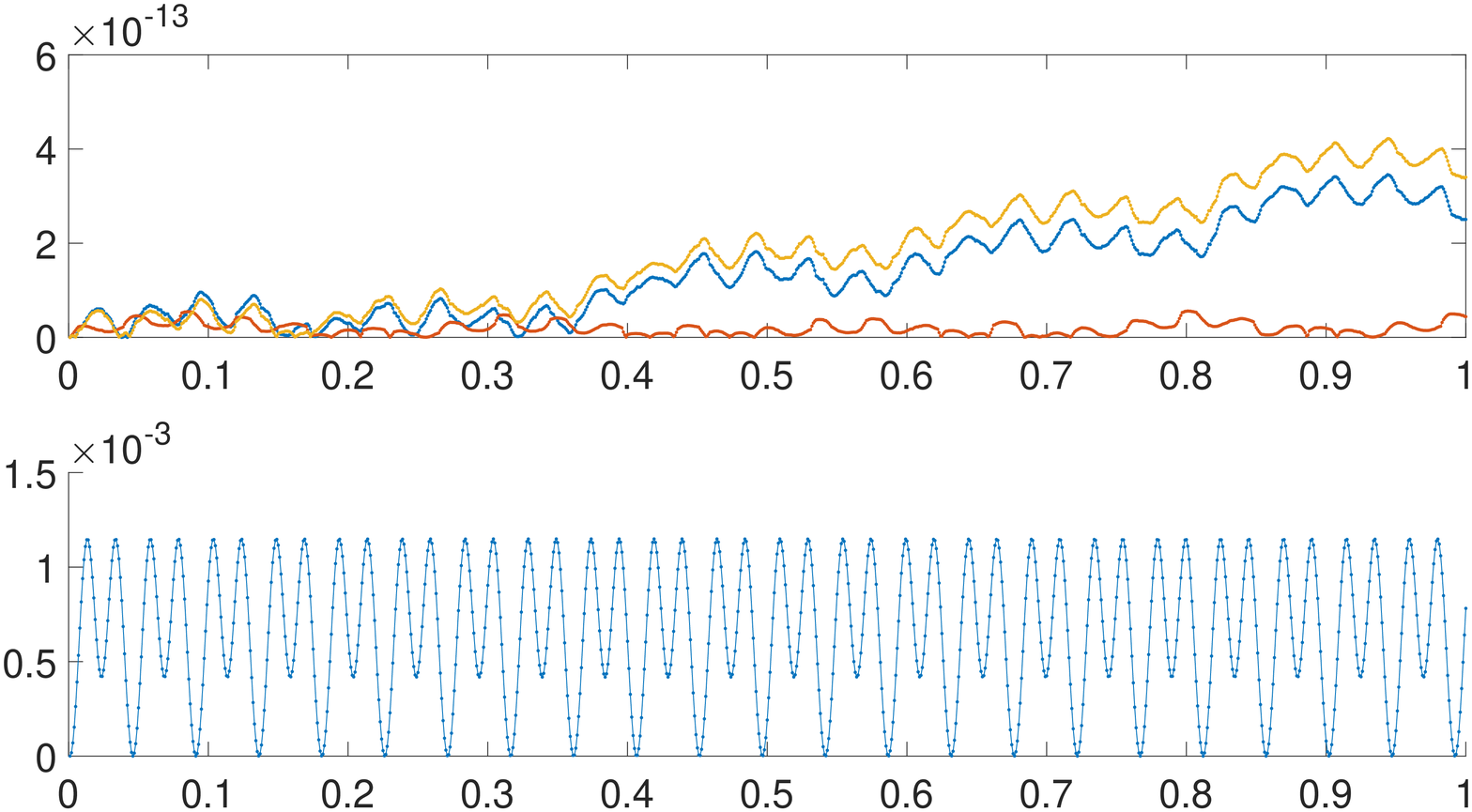}};
\centering
    \node[below=of img, node distance=0cm, yshift=4.1cm] {Hamiltonian variation in time for $w_0^2$};
  \node[below=of img, node distance=0cm, yshift=7.2cm] {Momentum variation in time  for $w_0^2$};
 \end{tikzpicture}
\end{minipage}
\caption{Momentum and Hamiltonian variation in time for the initial conditions $w_0^1,\Gamma^1$ and $w_0^2,\Gamma^2$, as defined above. The simulations have been carried with the hyperbolic midpoint method \eqref{IsoMP_min_sl2}, with time-step, respectively, $h=0.01$ and  $h=0.001$, total time $T=10$ and $T=1$ and tolerance for the Newton iteration $tol=10^{-13}$.}\label{fig:PV_ham_mom}
\end{figure}

We conclude showing in figure~\ref{fig:PV_ham_mom} the conservation properties for the hyperbolic midpoint method \eqref{IsoMP_min_sl2}, concerning the first integrals \eqref{eq:PV_mom} and \eqref{eq:PV_ham}. 

\clearpage 

\textbf{Acknowledgements}
The author would like to thank Klas Modin for the support and the enlightening discussions during the work on this paper.

\bibliographystyle{plainnat}
\bibliography{biblio}

 \end{document}